\pdfoutput=1   
\documentclass[pdflatex,sn-chicago]{sn-jnl}
\usepackage{amsmath,amssymb}
\usepackage{xcolor}
\usepackage{graphicx}
\usepackage{float}
\usepackage[section]{placeins}
\usepackage{array}
\usepackage{booktabs}
\usepackage{rotating}
\usepackage{tikz}
\usepackage{orcidlink}
\usetikzlibrary{calc,arrows.meta}
\hypersetup{hypertexnames=false}

\definecolor{paramcolor}{RGB}{0,70,140}

\title[Khayyam's Choices of Conics]{Khayyam's Choices of Conics}

\author*{\fnm{Amir} \sur{Asghari}\,\orcidlink{0000-0001-9658-2137}}\email{a.h.asghari@ljmu.ac.uk}
\affil{\orgname{Liverpool John Moores University}, \city{Liverpool}, \country{United Kingdom}}

\abstract{In the late eleventh century, Omar Khayyam gave a systematic classification of, and solutions to, cubic equations by means of intersecting conic sections. His work is often retold in modern Cartesian language, where each cubic is reduced to equations of certain curves whose intersection determines solutions in the Cartesian plane. In this paper, we retain the first part of the standard modern treatment of Khayyam: translating his verbal algebra into modern symbolic algebra. But instead of interpreting the resulting equations as curves in a Cartesian plane, we interpret them as Khayyam himself did, as relations between positive geometrical magnitudes represented by conics in the Apollonian sense. On this basis, and using a purpose-designed colour scheme, we reconstruct, in a compact treatment, his solutions to all the cubic equations treated by Khayyam. The reconstruction brings to light a hidden feature of the treatise and draws a sharp line between the geometrization of algebra as practised by Khayyam, on the one hand, and the geometrization of algebra developed by Fermat and Descartes, on the other.}

\keywords{Omar Khayyam, cubic equations, conic sections, analytic geometry, Apollonius, history of algebra}

\begin{document}
\maketitle

\definecolor{khblue}{RGB}{48,115,172}
\definecolor{khgreen}{RGB}{73,140,92}
\definecolor{khorange}{RGB}{222,139,45}
\definecolor{khred}{RGB}{185,55,55}
\tikzset{
  chosen/.style={black,line width=0.9pt},                                 
  unchosen/.style={black!70,dash pattern=on 3pt off 2pt,line width=0.7pt},
  localaxis/.style={black!60,line width=0.5pt},                           
  band/.style={line width=3.2pt,line cap=butt},
  xseg/.style={band,khred},      
  bseg/.style={band,khblue},     
  lseg/.style={band,khgreen},    
  aseg/.style={band,khorange},   
  yseg/.style={black,dash pattern=on 2.5pt off 2pt,line width=0.9pt},     
  solpoint/.style={circle,fill=black,inner sep=1.2pt},
  olab/.style={font=\small}
}
\newcommand{\localcross}[2]{%
  \draw[localaxis] (#1-0.18,#2)--(#1+0.18,#2);
  \draw[localaxis] (#1,#2-0.18)--(#1,#2+0.18);
  \node[olab,below left] at (#1,#2) {$O$};
}

\section{Introduction}

Written in the late eleventh century, Omar Khayyam's \emph{Treatise on Algebra} presents a systematic treatment of cubic equations by means of intersecting conic sections. Modern presentations frequently reformulate these constructions in Cartesian coordinates. Although such reformulations are mathematically legitimate and pedagogically convenient, they risk obscuring the geometric framework within which Khayyam himself worked. Khayyam begins his treatise with a warning: his work can be understood only by someone who has mastered Euclid's \emph{Elements} and Apollonius' \emph{Conics}. This paper accepts a similar constraint. It does not attempt to reconstruct the whole of Greek conic theory, nor to translate Khayyam prematurely into analytic geometry. Instead, it isolates the minimal geometric structures Khayyam repeatedly employs. It is within this framework, and with particular attention to Khayyam's choices of conics in each case, that the paper draws a fine line between Khayyam's geometry and the later Cartesian framework.

\section{Conics as geometric means}

Central to Khayyam's work is the geometric mean. Given two segments $a$ and $b$, a segment $x$ is their geometric mean if
\[
x^2=ab.
\]

Let $AB$ be a fixed segment and $H$ a variable point on it. At $H$ erect a perpendicular line, and choose a point $C$ on that perpendicular so that
\[
CH^2=AH\cdot HB.
\]
As $H$ varies along $AB$, the point $C$ traces a semicircle with diameter $AB$. In Euclidean language, the semicircle is the locus of points whose perpendicular to the diameter yields a segment that is the geometric mean of the two parts into which the diameter is divided.

If instead $H$ is taken on the continuation of $AB$ beyond $B$, and one requires
\[
CH^2=AB\cdot HB,
\]
then, as $H$ varies, the point $C$ traces a parabola.

If again $H$ is taken on the continuation of $AB$ beyond $B$, but one imposes
\[
CH^2=HA\cdot HB,
\]
the locus of $C$ is a hyperbola.

These definitions of the parabola and diameter-based hyperbola belong to Book I of Apollonius' \emph{Conics} \citep{ApolloniusHeath1896}. A complementary description, drawn from Book II, will be important later. For a hyperbola presented relative to its asymptotes, the rectangle formed by the perpendicular distances from a point on the curve to the two asymptotes is constant.

These constructions---circle, parabola, and hyperbola---form the geometric vocabulary within which Khayyam operates.

\section{Reading Khayyam}

Khayyam classifies cubic equations into thirteen species according to the arrangement of positive terms. He works exclusively with positive magnitudes interpreted as volumes; an equation equating a positive magnitude to zero is therefore excluded. This positivity constraint governs both classification and construction. He presents the cubic cases in two groups: first the species with three terms, then those with four. The latter divide further into two subgroups: cases with three terms on one side of the equation and one on the other, and cases with two terms on each side.

To see both what is gained and what can be lost when Khayyam's verbal algebra is translated into modern symbolic algebra, we begin by reading his first species closely. Its equality involves three volumetric terms---the sum of two is equal to a third---and hence belongs to the trinomial group. Unless otherwise stated, quotations from Khayyam, and the page numbers attached to them, follow the English translation of Rashed and Vahabzadeh \citep{RashedVahabzadeh2000}.

Khayyam announces the species as
\begin{quote}
\emph{``A cube and sides are equal to a number.''}
\end{quote}
We shall write this with all terms represented in the same dimension, here as volumes,
\[
x^3+b^2x=c^3.
\]
This preservation of dimensional type is what is meant here by \emph{homogeneity}. Our notation is deliberately geometrical: $x,b,c$ denote lengths, so that $x^3$ and $c^3$ are cubes and $b^2x$ is a solid. It differs superficially from the usual modern notation $x^3+bx=c$, but keeps the dimensional distinctions explicit in Khayyam's argument.

Khayyam begins:
\begin{quote}
\emph{``We suppose that $AB$ is the side of a square equal to the number of the roots.''}
\end{quote}

\begin{figure}[htbp]
\centering
\includegraphics[width=0.72\linewidth,height=0.54\textheight,keepaspectratio]{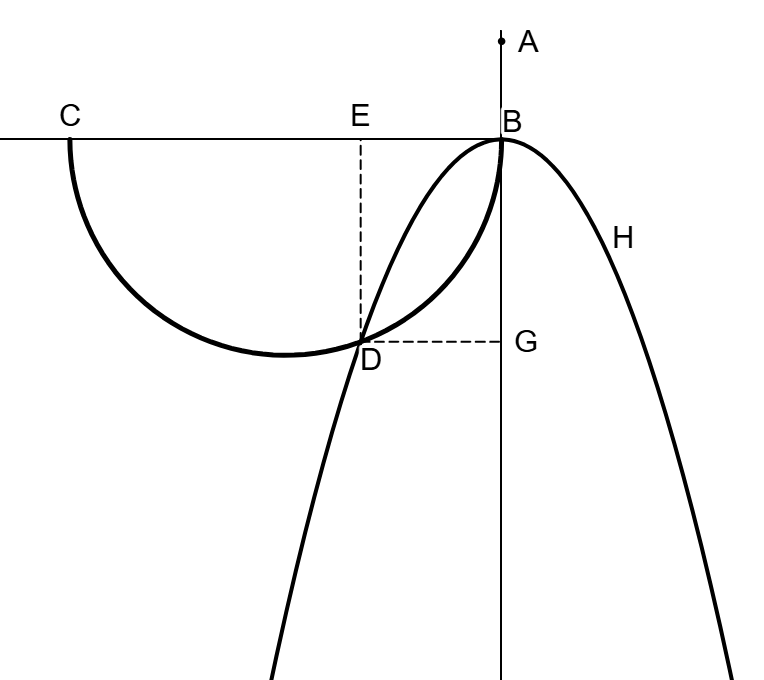}
\caption{Rashed and Vahabzadeh's reproduction of Khayyam's diagram accompanying their translation of the first species; Latin letters replace the letters of the Arabic original \citep[p.~131]{RashedVahabzadeh2000}}
\label{fig:rashed-reading}
\end{figure}

Thus, referring to Fig.~\ref{fig:rashed-reading}, in our notation
\[
AB=b,
\]
so that the square on $AB$ has area $b^2$, the coefficient of the ``sides'' in
\[
x^3+b^2x=c^3.
\]

He continues:
\begin{quote}
\emph{``And we construct a solid whose base is equal to the square of $AB$, whose height is equal to $BC$ and which is equal to the given number.''}
\end{quote}
The ``given number'' is the constant cube $c^3$. Hence
\[
AB^2\cdot BC=c^3.
\]
Since $AB=b$, we write
\[
BC=l,\qquad b^2l=c^3,
\]
and the cubic becomes
\[
x^3+b^2x=b^2l.
\]
Khayyam has therefore represented the constant solid on the same square base as the term containing the unknown.

He next constructs the first conic:
\begin{quote}
\emph{``And we produce $AB$ in a straight line to $G$. And we construct a parabola whose vertex is the point $B$, whose axis is $BG$ and whose erect side is $AB$.''}
\end{quote}
He then constructs the second:
\begin{quote}
\emph{``And we construct on $BC$ a semicircle. Then it will necessarily intersect the section: let them intersect at $D$.''}
\end{quote}
Finally,
\begin{quote}
\emph{``And we draw from $D$ \ldots\ two perpendiculars $DG,DE$ to $BG,BC$.''}
\end{quote}

At this point we introduce the two unknown lengths used in our modern retelling. The segment that will turn out to solve the cubic is
\[
BE=DG=x,
\]
while the companion segment is
\[
ED=BG=y.
\]
Since $BC=l$,
\[
EC=BC-BE=l-x.
\]

From the defining property of the parabola, whose erect side is $AB=b$,
\[
DG^2=BG\cdot AB,
\]
and hence
\[
x^2=by.
\]
From the semicircle on diameter $BC$,
\[
ED^2=BE\cdot EC,
\]
and therefore
\[
y^2=x(l-x).
\]
These equations record relations among Khayyam's constructed segments. The symbols $x,y,l$ and $b$ are simply names for positive lengths occurring in Khayyam's construction; so too are expressions such as $l-x$, since the construction itself ensures $x<l$. The same understanding will apply to analogous expressions in the remaining cases.

Khayyam then proves that the four segments $AB,BE,ED,EC$ are continuously proportional. In our notation,
\[
\frac{b}{x}=\frac{x}{y}=\frac{y}{l-x},
\]
or equivalently
\[
\frac{x}{b}=\frac{y}{x}=\frac{l-x}{y}.
\]
Consequently,
\[
\frac{x^2}{b^2}=\frac{l-x}{x},
\]
so
\[
x^3=b^2(l-x),
\]
and therefore
\[
x^3+b^2x=b^2l=c^3.
\]
Thus the segment $BE=x$ determined by the intersection of the parabola and semicircle is the required side of the cube.

Up to this point, our modernization differs little in substance from that of Rashed and Vahabzadeh. The difference in appearance comes mainly from notation. They begin with
\[
x^3+bx=c,
\]
so that $AB^2=b$ and $BC=c/b$, whereas we write
\[
x^3+b^2x=c^3=b^2l,
\]
taking $AB=b$ and $BC=l$. Their proportion
\[
\frac{b}{x^2}=\frac{x}{c/b-x}
\]
is therefore precisely our
\[
\frac{b^2}{x^2}=\frac{x}{l-x}.
\]
Likewise, their relations $BE^2=AB\cdot ED$ and $ED^2=BE\cdot EC$ become, under $BE=x$ and $ED=y$,
\[
x^2=by,\qquad y^2=x(l-x).
\]
Thus, up to notation, we have followed the same modernization.

Rashed and Vahabzadeh then explicitly discuss the alternatives available to Khayyam. They write:
\begin{quote}
``The \emph{choice} of the semicircle and the parabola is not the only choice possible: the parabola whose equation is $y=x^2+b$ could be considered, as could the hyperbola whose equation is $y=c/x$, by rewriting Equation 13 in the form $x(x^2+b)=c$; but this goes against the law of homogeneity.'' \emph{(italics ours)} \citep[p.~38]{RashedVahabzadeh2000}
\end{quote}
An alternative, however, does not require a separate algebraic rearrangement. It is already present in the continuous proportion itself. From
\[
\frac{x}{b}=\frac{y}{x}=\frac{l-x}{y},
\]
the three pairwise equalities yield
\[
x^2=by,\qquad y^2=x(l-x),\qquad x(y+b)=bl.
\]
Khayyam uses the first two relations, corresponding to the parabola and the circle, but any other pair is algebraically available---for example, the parabola and the hyperbola, as suggested by Rashed and Vahabzadeh.

No matter how one arrives at such an alternative pair, however, placing it on the same footing as the conics actually chosen by Khayyam already reflects a \emph{trained Cartesian habit}: treating algebraically available equations as curves equally available in one and the same plane. Our aim is instead to bracket that habit and ask what distinguishes Khayyam's choices within his own geometrical framework. Within that framework, homogeneity is unquestionably fundamental to Khayyam's algebra, but, as we shall see, it is not the reason for his choices of conics. First, however, let us see transparently how Khayyam gets to the conics in all thirteen cases.

\section{Khayyam's schemes}

The preceding section shows the characteristic direction of Khayyam's demonstrations. He introduces the conics one by one, identifies their point of intersection, derives the ratios satisfied by the segments determined by that point, and finally proves that one of those segments satisfies the original volumetric equation. The procedure establishes that the construction works, but it can leave a modern reader wondering how, starting from the cubic equation, one could arrive at the conics in the first place. This section reverses the direction of presentation and makes that derivation transparent. To avoid repeating the example just considered, we use Case~(7), which will also play an important role later.

Case~(7) is
\[
x^3+ax^2+b^2x=b^2l.
\]
Rearranging without leaving Khayyam's homogeneous setting gives
\[
x^2(x+a)=b^2(l-x),
\]
and hence
\[
\frac{x^2}{b^2}=\frac{l-x}{x+a}.
\]
Following a well-known Greek mean-proportional scheme, used for example in the duplication of the cube (for a historical account of that problem, see \citealp{Heath1921}), introduce an auxiliary positive segment $y$ and split this equality of ratios into the continuous proportion
\[
\frac{x}{b}=\frac{y}{x+a}=\frac{l-x}{y}.
\]
The three pairwise equalities yield three conic relations. From the first two terms,
\[
by=x(x+a).
\]
From the last two,
\[
y^2=(l-x)(x+a).
\]
And from the first and the last,
\[
xy=b(l-x),
\]
which is equivalently
\[
x(y+b)=bl.
\]
Thus the single proportional scheme produces
\[
\boxed{\;y^2=(l-x)(x+a),\qquad x(y+b)=bl,\qquad by=x(x+a).\;}
\]
Khayyam uses the first two relations: a circle and an asymptotic hyperbola. The third, a parabola, is algebraically available but is not used in his construction.

Case~(12) calls for a further explanation. In this species, the same symbolic cubic has to be read in two different positive forms according to the order of the given magnitudes $l$ and $a$:
\[
x^3+b^2x=ax^2+b^2l.
\]
When $l<a$, the relevant differences are positive in the form
\[
ax^2-x^3=b^2(x-l),
\]
whence
\[
x^2(a-x)=b^2(x-l),\qquad
\frac{x^2}{b^2}=\frac{x-l}{a-x},
\]
and the corresponding proportion
\[
\frac{x}{b}=\frac{y}{a-x}=\frac{x-l}{y}
\]
yields
\[
y^2=(x-l)(a-x),\qquad x(b-y)=bl,\qquad by=x(a-x).
\]
When $l>a$, the positive reading is instead
\[
x^3-ax^2=b^2(l-x),
\]
so that
\[
x^2(x-a)=b^2(l-x),\qquad
\frac{x^2}{b^2}=\frac{l-x}{x-a}.
\]
Introducing $y$ therefore gives
\[
\frac{x}{b}=\frac{y}{x-a}=\frac{l-x}{y},
\]
and hence
\[
y^2=(l-x)(x-a),\qquad x(y+b)=bl,\qquad by=x(x-a).
\]
The two sets of conic equations arise because Khayyam keeps every magnitude positive and chooses the difference that the geometrical order permits. When $l=a$, the cubic is satisfied directly by $x=a$.

The same procedure applies throughout the thirteen irreducible cubic species. Presenting the cases in this compact form separates the derivation of the relations from the later question of how the corresponding conics are constructed and placed. Khayyam sometimes subdivides a species further according to the relative magnitudes of the given data. In most such subdivisions the relations defining the chosen conics remain the same, while their relative placement and intersection behaviour change. Case~(12) is exceptional in the present list: positivity itself changes the form of the conic relations according as $l<a$ or $l>a$. The thirteen schemes are as follows.

\begingroup
\abovedisplayskip=4pt
\belowdisplayskip=3pt
\abovedisplayshortskip=2pt
\belowdisplayshortskip=2pt

\[
\begin{aligned}
(1)\ &x^3+b^2x=b^2l\\[-2pt]
&x^2=by,\quad y^2=x(l-x)\\[-2pt]
&x(y+b)=bl
\end{aligned}
\]
\[
\begin{aligned}
(2)\ &x^3+b^2l=b^2x\\[-2pt]
&x^2=by,\quad y^2=x(x-l)\\[-2pt]
&x(b-y)=bl
\end{aligned}
\]
\[
\begin{aligned}
(3)\ &x^3=b^2x+b^2l\\[-2pt]
&x^2=by,\quad y^2=x(l+x)\\[-2pt]
&x(y-b)=bl
\end{aligned}
\]
\[
\begin{aligned}
(4)\ &x^3+ax^2=c^3\\[-2pt]
&xy=c^2,\quad y^2=c(x+a)\\[-2pt]
&x(x+a)=cy
\end{aligned}
\]
\[
\begin{aligned}
(5)\ &x^3+c^3=ax^2\\[-2pt]
&xy=c^2,\quad y^2=c(a-x)\\[-2pt]
&cy=x(a-x)
\end{aligned}
\]
\[
\begin{aligned}
(6)\ &x^3=ax^2+c^3\\[-2pt]
&xy=c^2,\quad y^2=c(x-a)\\[-2pt]
&cy=x(x-a)
\end{aligned}
\]
\[
\begin{aligned}
(7)\ &x^3+ax^2+b^2x=b^2l\\[-2pt]
&y^2=(l-x)(x+a),\quad x(y+b)=bl\\[-2pt]
&by=x(x+a)
\end{aligned}
\]
\[
\begin{aligned}
(8)\ &x^3+ax^2+b^2l=b^2x\\[-2pt]
&y^2=(x-l)(x+a),\quad x(b-y)=bl\\[-2pt]
&by=x(x+a)
\end{aligned}
\]
\[
\begin{aligned}
(9)\ &x^3+b^2x+b^2l=ax^2\\[-2pt]
&y^2=(x+l)(a-x),\quad x(y-b)=bl\\[-2pt]
&by=x(a-x)
\end{aligned}
\]
\[
\begin{aligned}
(10)\ &x^3=ax^2+b^2x+b^2l\\[-2pt]
&y^2=(x+l)(x-a),\quad x(y-b)=bl\\[-2pt]
&by=x(x-a)
\end{aligned}
\]
\[
\begin{aligned}
(11)\ &x^3+ax^2=b^2x+b^2l\\[-2pt]
&y^2=(x+l)(x+a),\quad x(y-b)=bl\\[-2pt]
&by=x(x+a)
\end{aligned}
\]
\[
\begin{aligned}
(12)\ &x^3+b^2x=ax^2+b^2l\\[-2pt]
&l<a:\quad y^2=(x-l)(a-x),\quad x(b-y)=bl,\quad by=x(a-x)\\[-2pt]
&l>a:\quad y^2=(l-x)(x-a),\quad x(y+b)=bl,\quad by=x(x-a)\\[-2pt]
&l=a:\quad x=a\quad\text{directly.}
\end{aligned}
\]
\[
\begin{aligned}
(13)\ &x^3+b^2l=ax^2+b^2x\\[-2pt]
&y^2=(x-l)(x-a),\quad x(b-y)=bl\\[-2pt]
&by=x(x-a)\\[-2pt]
&l=a:\quad x=a\quad\text{directly.}
\end{aligned}
\]
\endgroup

\section{Local coordinates without global coordinates}

It is tempting, from a modern perspective, to read these constructions in the Cartesian plane: fix an origin, fix perpendicular axes, and then consider the curves defined by the given equations. This translation is mathematically legitimate. But it imports a framework that Khayyam does not use. The Cartesian framework makes position, sign, and shift appear automatic, whereas in Khayyam's geometry the orientation and relocation of the defining segments must be handled directly within the construction.

To see the difference sharply, consider the following two parabolas written in Cartesian form:
\[
y=x^2,\qquad y=-(x-a)^2+b.
\]
In Cartesian language these are distinct equations; one of them ``opens downward'' because of the minus sign. Their graphs are positioned differently in the plane, and their algebraic forms encode that difference.

In the Apollonian framework, however, this distinction is largely beside the point. The appearance of a negative sign in a Cartesian equation is not something that exists for Khayyam at all; it is an artifact of global coordinate conventions. The curve is the curve. Its Apollonian description does not change when we translate it, rotate it, or relabel which direction is called ``up.'' What matters is not position in a fixed plane, but the relation that defines the curve.

Once we agree that the diameter is perpendicular to the tangent at the vertex, the only datum distinguishing one parabola from another is the parameter $p$ in
\[
x^2=py \qquad (\text{equivalently } y^2=px).
\]
The number $p$ represents the length of a segment fixing how open the curve is. In this sense, $p$ plays a role analogous to the radius of a circle: changing $p$ changes the curve itself, whereas translating the curve in the plane does not.

In Khayyam's work, each conic is introduced through a local construction that comes with the conic itself. One first specifies the lines that play the role of a coordinate cross for that conic---its diameter and tangent at the vertex, or, in the asymptotic case, its two asymptotes---and only then measures distances relative to those lines. The ``coordinates'' are not signed distances from fixed axes in a global plane, but lengths of constructed segments.

Because Khayyam works in local coordinates, when two curves visually meet, it does not automatically follow that the intersection determines the solution of the relevant equation; and, conversely, when they do not meet in a given placement, it does not follow that the equation has no solution. Everything depends on how meaningfully the curves have been placed within the local configuration. Figures~\ref{fig:local-meet}--\ref{fig:global-meet} illustrate the point with the relations $y=x^2$ and $x^2+y^2=1$. For the circle no analogous coordinate cross is required: its relevant local data are supplied by its diameter.

\begin{figure}[htbp]
\centering
\begin{tikzpicture}[scale=1.18]
  \draw[localaxis] (-0.62,0)--(0.62,0);
  \draw[localaxis] (0,-0.24)--(0,1.22);
  \draw[chosen,domain=-1.38:1.38,samples=90,smooth,variable=\t] plot ({\t},{0.70*\t*\t});
  \coordinate (P) at (0.80,0.448);
  \coordinate (C) at (1.36,0.91);
  \draw[chosen] (C) circle (0.726);
  \node[solpoint] at (P) {};
\end{tikzpicture}
\caption{The conics meet visually, but the intersection cannot be read as a common solution: the parabola and circle have not been placed so that the relevant constructed lengths coincide}
\label{fig:local-meet}
\end{figure}

\begin{figure}[htbp]
\centering
\begin{tikzpicture}[scale=1.05]
  \begin{scope}[shift={(-1.55,-0.55)},rotate=8]
    \draw[localaxis] (-0.58,0)--(0.58,0);
    \draw[localaxis] (0,-0.22)--(0,1.05);
    \draw[chosen,domain=-1.35:1.35,samples=80,smooth,variable=\t] plot ({\t},{0.72*\t*\t});
  \end{scope}
  \begin{scope}[shift={(1.65,1.05)},rotate=-18]
    \draw[chosen] (0,0) circle (1.05);
  \end{scope}
\end{tikzpicture}
\caption{The same two local relations can be placed so that the conics do not meet at all; this does not imply that the corresponding equation has no solution}
\label{fig:local-no-meet}
\end{figure}

\begin{figure}[htbp]
\centering
\begin{tikzpicture}[scale=2.0]
  \draw[->,localaxis] (-1.25,0)--(1.35,0) node[right] {$x$};
  \draw[->,localaxis] (0,-1.25)--(0,1.45) node[above] {$y$};
  \node[olab,below left] at (0,0) {$O$};
  \draw[chosen,domain=-1.08:1.08,samples=100,smooth,variable=\t] plot ({\t},{\t*\t});
  \draw[chosen] (0,0) circle (1);
  \node[solpoint] at (0.78615,0.61803) {};
  \draw[yseg] (0.78615,0)--(0.78615,0.61803);
  \draw[xseg] (0,0)--(0.78615,0);
  \node[font=\small,anchor=west] at (0.93,1.18) {$y=x^2$};
  \node[font=\small,anchor=east] at (-1.03,0.72) {$x^2+y^2=1$};
\end{tikzpicture}
\caption{In a global Cartesian frame the two equations share the same $x$ and $y$ by definition; their intersection therefore carries a common pair of coordinates}
\label{fig:global-meet}
\end{figure}
\FloatBarrier

The contrast is not that Khayyam lacked coordinate-like measurements. Rather, the common frame itself had to be constructed. A global Cartesian plane supplies that compatibility in advance; Khayyam's geometry does not.

\section{Khayyam's choices of conics}

Equipped with this local framework, we can now return to Khayyam's choices of conics and ask why, in each case, he used the two he did rather than the third relation algebraically available to him.

In the colour figures in this section, the colours record constructed magnitudes: blue for $b$ and $c$, green for $l$, orange for $a$, red for the sought segment $x$, and a dashed black segment for its companion $y$. This use of colour also follows the spirit of the exposition by Kent and Muraki, who replace many letter labels by colour in presenting one of Khayyam's cubic constructions \citep{KentMuraki2016}. Each figure shows two thin reference lines: the diameter line carrying the given segments, and its perpendicular at $O$.

The label $O$ should not be read as the origin of a global coordinate plane. It marks the local reference point from which the red segment $x$ is read and, when the visible cross belongs to a non-circular conic, it identifies the defining local frame of that conic. For a parabola or a diameter-based hyperbola, the two lines of the cross are its diameter and the tangent at the vertex; for an asymptotic hyperbola, they are its asymptotes. A circle needs no analogous cross: its relevant local data are carried by its diameter, and $O$ is simply the chosen point on that diameter from which the required segment is measured.

The two thin lines just described are in each figure lines of these local crosses: the perpendicular at $O$ is the axis of the parabola in Cases~(1)--(3) and an asymptote of the asymptotic hyperbola in the remaining cases, while the diameter line serves as tangent at the vertex, second asymptote, or diameter according to the case. No cross is completed beyond these two lines. The reader can nevertheless recover the local cross belonging to each of the other non-circular conics from the construction. This is important when reading the algebra back from the picture: each conic has its own diameter-and-tangent or pair of asymptotes, even when those lines are not drawn to full length. In particular, the blue segment $b$ often records the displacement to an asymptote of the asymptotic hyperbola. If the horizontal line through the far end of that blue segment were drawn, it would in many of the figures be one of the hyperbola's asymptotes. Imagining that omitted line makes expressions such as $y+b$ or $y-b$ directly visible. We suppress these completed crosses and extended asymptotes in order to keep the figures to the essential geometry and to avoid cluttering the constructions.

Where Khayyam distinguishes several configurations within a species, the reconstruction shows a configuration yielding the required intersection; the captions indicate this where necessary. Except in Case~(12), such variations alter the relative placement and intersection behaviour of the same conic relations rather than the relations themselves.

Let us begin with Case~(1) (Fig.~\ref{fig:case1}). The two relations used by Khayyam are
\[
x^2=by,\qquad y^2=x(l-x),
\]
while the remaining relation is
\[
x(y+b)=bl.
\]
The chosen parabola supplies the guiding construction. The semicircle can be fitted to it so that the same red segment $x$ and dashed segment $y$ are read directly in both. If, instead, the rectangular hyperbola were used, the dashed segment $y$ would have to be extended by the already known blue segment $b$ (that is, $y+b$). The alternative is possible, but it asks the demonstration to move between the original dashed segment and its extension by $b$.

  \begin{figure}[htbp]
  \centering
\begin{tikzpicture}[x=1.422cm,y=1.422cm]
\clip (-6.2,-4.6) rectangle (0.55,1.25);
\draw[localaxis] (0,0) -- (-4.4,0);
\draw[localaxis] (0,0) -- (0,-4.6);
\draw[unchosen,domain=-0.3267:1.9021,samples=161,smooth,variable=\t] plot ({-(exp(\t))},{-(-1+4.4/exp(\t))});
\draw[chosen,domain=0:180,samples=91,smooth,variable=\t] plot ({-(2.2+2.2*cos(\t))},{-(2.2*sin(\t))});
\draw[chosen,domain=-2.3083:2.3083,samples=121,smooth,variable=\t] plot ({-(\t)},{-(\t*\t/1)});
\draw[lseg,line width=3.2pt,yshift=-1.85pt] (0,0) -- (-4.4,0);
\draw[bseg] (0,0) -- (0,1);
\draw[xseg,yshift=1.85pt] (0,0) -- (-1.4364,0);
\draw[yseg] (-1.4364,0) -- (-1.4364,-2.0632);
\node[solpoint] at (0,0) {};
\node[olab,anchor=north,xshift=7pt,yshift=-1.85pt] at (0,0) {$O$};
\end{tikzpicture}
  \caption{Case (1): the two conics chosen by Khayyam and the remaining rectangular hyperbola}\label{fig:case1}
  \end{figure}%

Cases~(2) and (3) exhibit the same kind of adjustment (Figs.~\ref{fig:case2} and \ref{fig:case3}).

  \begin{figure}[htbp]
  \centering
\begin{tikzpicture}[x=3.652cm,y=3.652cm]
\clip (-1.27,-1.17) rectangle (0.56,1.13);
\draw[localaxis] (0,0) -- (-0.4289,0);
\draw[localaxis] (0,0) -- (0,-1.17);
\draw[unchosen,domain=-2.0168:0.571,samples=161,smooth,variable=\t] plot ({-(exp(\t))},{-(1-0.35/exp(\t))});
\draw[chosen,domain=-2.5678:2.5678,samples=141,smooth,variable=\t] plot ({-(0.175+0.175*cosh(\t))},{-(0.175*sinh(\t))});
\draw[chosen,domain=-1.3423:1.3423,samples=121,smooth,variable=\t] plot ({-(\t)},{-(\t*\t/1)});
\draw[lseg,line width=3.2pt,yshift=-1.85pt] (0,0) -- (-0.35,0);
\draw[bseg] (0,0) -- (0,-1);
\draw[xseg,yshift=1.85pt] (0,0) -- (-0.4289,0);
\draw[yseg] (-0.4289,0) -- (-0.4289,-0.184);
\node[solpoint] at (0,0) {};
\node[olab,anchor=north,xshift=7pt,yshift=-1.85pt] at (0,0) {$O$};
\end{tikzpicture}
  \caption{Case (2): one solvable configuration of Khayyam's chosen conics and the remaining asymptotic hyperbola}\label{fig:case2}
  \end{figure}%

  \begin{figure}[htbp]
  \centering
\begin{tikzpicture}[x=1.235cm,y=1.235cm]
\clip (-1.4,-4) rectangle (4.2,2.8);
\draw[localaxis] (-1.15,0) -- (1.3588,0);
\draw[localaxis] (0,0) -- (0,-4);
\draw[unchosen,domain=-1.3189:1.5476,samples=161,smooth,variable=\t] plot ({exp(\t)},{-(1+1.15/exp(\t))});
\draw[chosen,domain=-2.3864:2.7096,samples=141,smooth,variable=\t] plot ({-0.575+0.575*cosh(\t)},{-(0.575*sinh(\t))});
\draw[chosen,domain=-2.1713:2.1713,samples=121,smooth,variable=\t] plot ({\t},{-(\t*\t/1)});
\draw[lseg,line width=3.2pt,yshift=-1.85pt] (0,0) -- (-1.15,0);
\draw[bseg] (0,0) -- (0,-1);
\draw[xseg,yshift=1.85pt] (0,0) -- (1.3588,0);
\draw[yseg] (1.3588,0) -- (1.3588,-1.8463);
\node[solpoint] at (0,0) {};
\node[olab,anchor=north,xshift=-7pt,yshift=-5.05pt] at (0,0) {$O$};
\end{tikzpicture}
  \caption{Case (3): Khayyam's chosen conics and the remaining asymptotic hyperbola}\label{fig:case3}
  \end{figure}%

In Cases~(1)--(3), the departure required by the remaining conic is comparatively modest: the dashed segment is extended or shortened by the already available blue segment $b$. From Case~(4) onwards (Figs.~\ref{fig:case4}--\ref{fig:case6}), the contrast becomes stronger. The conics Khayyam chooses may require segments such as $x\pm a$, $y\pm b$, or $l\pm x$, but these are assembled directly from segments already present and \emph{known in position}---the phrase Khayyam uses throughout the treatise whenever he records the data for drawing a conic or introduces a segment or a point, including the segment that solves the problem \citep{RashedVahabzadeh2000}. The remaining parabola stands differently: its vertex lies at $a/2$ from $O$, displaced by $a^2/(4c)$ in Cases~(4)--(6) and by $a^2/(4b)$ in Cases~(7)--(13). The point is not that these segments were hard to construct---bisecting a segment and finding a third proportional were entirely within Khayyam's reach. The point is that, until such a construction is carried out, the cross of the remaining conic is not known in position: it stands at no point already marked in the figure. Algebraically, completing the square makes the placement immediate; geometrically, it is a further act of construction.

  \begin{figure}[htbp]
  \centering
\begin{tikzpicture}[x=1.508cm,y=1.508cm]
\clip (-3.9,-1.17) rectangle (2.2,4.4);
\draw[localaxis] (2,0) -- (-0.618,0);
\draw[localaxis] (0,0) -- (0,4.4);
\draw[unchosen,domain=-2.4:4.1,samples=161,smooth,variable=\t] plot ({-(\t)},{\t*(\t+2)/1});
\draw[chosen,domain=-1.5892:1.4816,samples=161,smooth,variable=\t] plot ({-(exp(\t))},{0+1/exp(\t)});
\draw[chosen,domain=0:2.49,samples=121,smooth,variable=\t] plot ({-(-2+\t*\t/1)},{\t});
\draw[aseg,line width=3.2pt,yshift=-1.85pt] (0,0) -- (2,0);
\draw[bseg] (0,0) -- (0,1);
\draw[xseg,yshift=1.85pt] (0,0) -- (-0.618,0);
\draw[yseg] (-0.618,0) -- (-0.618,1.618);
\node[solpoint] at (0,0) {};
\node[olab,anchor=north,xshift=10pt,yshift=-5.05pt] at (0,0) {$O$};
\end{tikzpicture}
  \caption{Case (4): Khayyam's chosen conics and the remaining parabola}\label{fig:case4}
  \end{figure}%

  \begin{figure}[htbp]
  \centering
\begin{tikzpicture}[x=1.032cm,y=1.032cm]
\clip (-2.5,-2) rectangle (6.8,5);
\draw[localaxis] (0,0) -- (3.3,0);
\draw[localaxis] (0,0) -- (0,5);
\draw[unchosen,domain=-2.7:7,samples=161,smooth,variable=\t] plot ({\t},{-\t*(\t-3.3)/1});
\draw[chosen,domain=-1.7047:1.9879,samples=161,smooth,variable=\t] plot ({exp(\t)},{0+1/exp(\t)});
\draw[chosen,domain=0:2.4698,samples=121,smooth,variable=\t] plot ({3.3-\t*\t/1},{\t});
\draw[aseg,line width=3.2pt,yshift=-1.85pt] (0,0) -- (3.3,0);
\draw[bseg] (0,0) -- (0,1);
\draw[xseg,yshift=1.85pt] (0,0) -- (0.6097,0);
\draw[yseg] (0.6097,0) -- (0.6097,1.6402);
\node[solpoint] at (0,0) {};
\node[olab,anchor=north,xshift=-10pt,yshift=-1.85pt] at (0,0) {$O$};
\end{tikzpicture}
  \caption{Case (5): one solvable configuration of Khayyam's chosen conics and the remaining parabola; the defining conic relations are unchanged in Khayyam's other parameter configurations}\label{fig:case5}
  \end{figure}%

  \begin{figure}[htbp]
  \centering
\begin{tikzpicture}[x=1.455cm,y=1.455cm]
\clip (-3.7,-1) rectangle (2.9,3.9);
\draw[localaxis] (0,0) -- (-1.733,0);
\draw[localaxis] (0,0) -- (0,3.9);
\draw[unchosen,domain=-3.1:3.9,samples=161,smooth,variable=\t] plot ({-(\t)},{\t*(\t-1.4)/1});
\draw[chosen,domain=-1.4816:1.4351,samples=161,smooth,variable=\t] plot ({-(exp(\t))},{0+1/exp(\t)});
\draw[chosen,domain=0:1.6125,samples=121,smooth,variable=\t] plot ({-(1.4+\t*\t/1)},{\t});
\draw[aseg,line width=3.2pt,yshift=-1.85pt] (0,0) -- (-1.4,0);
\draw[bseg] (0,0) -- (0,1);
\draw[xseg,yshift=1.85pt] (0,0) -- (-1.733,0);
\draw[yseg] (-1.733,0) -- (-1.733,0.577);
\node[solpoint] at (0,0) {};
\node[olab,anchor=north,xshift=7pt,yshift=-1.85pt] at (0,0) {$O$};
\end{tikzpicture}
  \caption{Case (6): Khayyam's chosen conics and the remaining parabola}\label{fig:case6}
  \end{figure}%

Case~(7) (Fig.~\ref{fig:case7}) makes the difference especially visible. The first relation is
\[
y^2=(l-x)(x+a).
\]
In the intended colour figure, the dashed perpendicular $y$ is the geometric mean of the two portions of the diameter represented by $l-x$ and $x+a$. The sought red segment $x$ is therefore read together with the already constructed orange segment $a$. The second chosen relation,
\[
x(y+b)=bl,
\]
uses the same red segment and extends the dashed segment by the already known blue segment $b$. Both adjustments are already present in the figure.

The remaining relation is
\[
x(x+a)=by.
\]
To construct its parabola in the same setting, one must first locate a new vertex. In modern symbols the required displacement becomes visible by completing the square,
\[
\left(x+\frac a2\right)^2=b\left(y+\frac{a^2}{4b}\right),
\]
and each of the steps it calls for was well within Khayyam's reach. Had he chosen this parabola, his own idiom suggests how the construction would have run: the midpoint of the segment $a$ is known in position; at that midpoint a perpendicular to the diameter is erected; and on that perpendicular a segment equal to $a^2/(4b)$ is constructed. Only now is the cross of the parabola known in position, and with it the parabola itself. The issue is therefore not the difficulty of these steps but their necessity: the crosses of the two chosen conics are known in position from the data already laid out, while the cross of the remaining conic becomes so only through further construction. This is the extra geometrical work that the global algebraic form hides. Figures~\ref{fig:case8}--\ref{fig:case13-gt} show the corresponding configurations for the remaining cases.

  \begin{figure}[htbp]
  \centering
\begin{tikzpicture}[x=1.745cm,y=1.745cm]
\clip (-2.8,-1.35) rectangle (2.7,3.2);
\draw[localaxis] (0.8,0) -- (-1.6,0);
\draw[localaxis] (0,0) -- (0,3.2);
\draw[unchosen,domain=-2.9:3,samples=161,smooth,variable=\t] plot ({-(\t)},{\t*(\t+0.8)/1});
\draw[chosen,domain=0:180,samples=91,smooth,variable=\t] plot ({-(0.4+1.2*cos(\t))},{1.2*sin(\t)});
\draw[chosen,domain=-1.0776:1.1939,samples=161,smooth,variable=\t] plot ({-(exp(\t))},{-1+1.6/exp(\t)});
\draw[lseg,line width=3.2pt,yshift=-1.85pt] (0,0) -- (-1.6,0);
\draw[aseg,line width=3.2pt,yshift=-1.85pt] (0,0) -- (0.8,0);
\draw[bseg] (0,0) -- (0,-1);
\draw[xseg,yshift=1.85pt] (0,0) -- (-0.7443,0);
\draw[yseg] (-0.7443,0) -- (-0.7443,1.1495);
\node[solpoint] at (0,0) {};
\node[olab,anchor=north,xshift=7pt,yshift=-5.05pt] at (0,0) {$O$};
\end{tikzpicture}
  \caption{Case (7): circle and asymptotic hyperbola chosen by Khayyam, with the remaining parabola}\label{fig:case7}
  \end{figure}%

  \begin{figure}[htbp]
  \centering
\begin{tikzpicture}[x=2.233cm,y=2.233cm]
\clip (-2.35,-0.35) rectangle (1.95,2.9);
\draw[localaxis] (0.44,0) -- (-0.637,0);
\draw[localaxis] (0,-0.35) -- (0,2.9);
\draw[unchosen,domain=-2.15:2.55,samples=161,smooth,variable=\t] plot ({-(\t)},{\t*(\t+0.44)/1});
\draw[chosen,domain=-1.4575:2.771,samples=141,smooth,variable=\t] plot ({-(-0.12+0.32*cosh(\t))},{0.32*sinh(\t)});
\draw[chosen,domain=-2.4849:1.0473,samples=161,smooth,variable=\t] plot ({-(exp(\t))},{1-0.2/exp(\t)});
\draw[lseg,line width=3.2pt,yshift=-1.85pt] (0,0) -- (-0.2,0);
\draw[aseg,line width=3.2pt,yshift=-1.85pt] (0,0) -- (0.44,0);
\draw[bseg] (0,0) -- (0,1);
\draw[xseg,yshift=1.85pt] (0,0) -- (-0.637,0);
\draw[yseg] (-0.637,0) -- (-0.637,0.686);
\node[solpoint] at (0,0) {};
\node[olab,anchor=north,xshift=7pt,yshift=-5.05pt] at (0,0) {$O$};
\end{tikzpicture}
  \caption{Case (8): one solvable configuration of Khayyam's chosen conics and the remaining parabola; Khayyam also considers tangency and the absence of a required further intersection}\label{fig:case8}
  \end{figure}%

  \begin{figure}[htbp]
  \centering
\begin{tikzpicture}[x=1.000cm,y=1.000cm]
\clip (-6.1,-2) rectangle (2.5,6.4);
\draw[localaxis] (2,0) -- (-4,0);
\draw[localaxis] (0,0) -- (0,6.4);
\draw[unchosen,domain=-2.7:6.3,samples=161,smooth,variable=\t] plot ({-(\t)},{-\t*(\t-4)/1});
\draw[chosen,domain=0:180,samples=91,smooth,variable=\t] plot ({-(1+3*cos(\t))},{3*sin(\t)});
\draw[chosen,domain=-1.0818:1.8871,samples=161,smooth,variable=\t] plot ({-(exp(\t))},{1+2/exp(\t)});
\draw[aseg,line width=3.2pt,yshift=-1.85pt] (0,0) -- (-4,0);
\draw[lseg,line width=3.2pt,yshift=-1.85pt] (0,0) -- (2,0);
\draw[bseg] (0,0) -- (0,1);
\draw[xseg,yshift=1.85pt] (0,0) -- (-1,0);
\draw[yseg] (-1,0) -- (-1,3);
\node[solpoint] at (0,0) {};
\node[olab,anchor=north,xshift=10pt,yshift=-5.05pt] at (0,0) {$O$};
\end{tikzpicture}
  \caption{Case (9): one solvable configuration of Khayyam's chosen conics and the remaining parabola}\label{fig:case9}
  \end{figure}%

  \begin{figure}[htbp]
  \centering
\begin{tikzpicture}[x=1.414cm,y=1.414cm]
\clip (-3.84,-1.24) rectangle (1.74,4.7);
\draw[localaxis] (0.8,0) -- (-2.1054,0);
\draw[localaxis] (0,0) -- (0,4.7);
\draw[unchosen,domain=-1.94:4.04,samples=161,smooth,variable=\t] plot ({-(\t)},{\t*(\t-1.45)/1});
\draw[chosen,domain=-1.1198:1.9021,samples=141,smooth,variable=\t] plot ({-(0.325+1.125*cosh(\t))},{1.125*sinh(\t)});
\draw[chosen,domain=-1.6582:1.4679,samples=161,smooth,variable=\t] plot ({-(exp(\t))},{1+0.8/exp(\t)});
\draw[aseg,line width=3.2pt,yshift=-1.85pt] (0,0) -- (-1.45,0);
\draw[lseg,line width=3.2pt,yshift=-1.85pt] (0,0) -- (0.8,0);
\draw[bseg] (0,0) -- (0,1);
\draw[xseg,yshift=1.85pt] (0,0) -- (-2.1054,0);
\draw[yseg] (-2.1054,0) -- (-2.1054,1.38);
\node[solpoint] at (0,0) {};
\node[olab,anchor=north,xshift=7pt,yshift=-5.05pt] at (0,0) {$O$};
\end{tikzpicture}
  \caption{Case (10): Khayyam's chosen conics and the remaining parabola}\label{fig:case10}
  \end{figure}%

  \begin{figure}[htbp]
  \centering
\begin{tikzpicture}[x=1.060cm,y=1.060cm]
\clip (-4.9,-1.6) rectangle (4.16,5.6);
\draw[localaxis] (2.3,0) -- (-0.6806,0);
\draw[localaxis] (0,0) -- (0,5.6);
\draw[unchosen,domain=-4.36:5.1,samples=161,smooth,variable=\t] plot ({-(\t)},{\t*(\t+2.3)/1});
\draw[chosen,domain=0:2.6958,samples=141,smooth,variable=\t] plot ({-(-1.5+0.8*cosh(\t))},{0.8*sinh(\t)});
\draw[chosen,domain=-1.9859:1.6864,samples=161,smooth,variable=\t] plot ({-(exp(\t))},{1+0.7/exp(\t)});
\draw[aseg,line width=4.8pt,yshift=-2.65pt] (0,0) -- (2.3,0);
\draw[bseg] (0,0) -- (0,1);
\draw[lseg,line width=2pt,yshift=-1.25pt] (0,0) -- (0.7,0);
\draw[xseg,yshift=1.85pt] (0,0) -- (-0.6806,0);
\draw[yseg] (-0.6806,0) -- (-0.6806,2.0285);
\node[solpoint] at (0,0) {};
\node[olab,anchor=north,xshift=12pt,yshift=-6.65pt] at (0,0) {$O$};
\end{tikzpicture}
  \caption{Case (11): one of Khayyam's configurations; changing the order of $l$ and $a$ changes the placement, not the defining conic relations}\label{fig:case11}
  \end{figure}%

  \begin{figure}[htbp]
  \centering
\begin{tikzpicture}[x=2.963cm,y=2.963cm]
\clip (-2.47,-1.49) rectangle (0.77,1.15);
\draw[localaxis] (0,0) -- (-1.75,0);
\draw[localaxis] (0,-1.49) -- (0,1.15);
\draw[unchosen,domain=-0.97:2.67,samples=161,smooth,variable=\t] plot ({-(\t)},{-\t*(\t-1.75)/1});
\draw[chosen,domain=0:180,samples=91,smooth,variable=\t] plot ({-(1.075+0.675*cos(\t))},{0.675*sin(\t)});
\draw[chosen,domain=-2.0116:1.0886,samples=161,smooth,variable=\t] plot ({-(exp(\t))},{1-0.4/exp(\t)});
\draw[aseg,line width=4.8pt,yshift=-2.65pt] (0,0) -- (-1.75,0);
\draw[bseg] (0,0) -- (0,1);
\draw[lseg,line width=2pt,yshift=-1.25pt] (0,0) -- (-0.4,0);
\draw[xseg,yshift=1.85pt] (0,0) -- (-1.1928,0);
\draw[yseg] (-1.1928,0) -- (-1.1928,0.6646);
\node[solpoint] at (0,0) {};
\node[olab,anchor=north,xshift=11pt,yshift=-1.85pt] at (0,0) {$O$};
\end{tikzpicture}
  \caption{Case (12), $l<a$: Khayyam's chosen conics and the remaining parabola}\label{fig:case12-lt}
  \end{figure}%

  \begin{figure}[htbp]
  \centering
\begin{tikzpicture}[x=1.231cm,y=1.231cm]
\clip (-5.6,-4.6) rectangle (2.2,1.6);
\draw[localaxis] (0,0) -- (-4,0);
\draw[localaxis] (0,0) -- (0,-4.6);
\draw[unchosen,domain=-2.4:5.8,samples=161,smooth,variable=\t] plot ({-(\t)},{-(\t*(\t-1)/1)});
\draw[chosen,domain=0:180,samples=91,smooth,variable=\t] plot ({-(2.5+1.5*cos(\t))},{-(1.5*sin(\t))});
\draw[chosen,domain=-0.422:1.8083,samples=161,smooth,variable=\t] plot ({-(exp(\t))},{-(-1+4/exp(\t))});
\draw[lseg,line width=4.8pt,yshift=-2.65pt] (0,0) -- (-4,0);
\draw[bseg] (0,0) -- (0,1);
\draw[aseg,line width=2pt,yshift=-1.25pt] (0,0) -- (-1,0);
\draw[xseg,yshift=1.85pt] (0,0) -- (-1.743,0);
\draw[yseg] (-1.743,0) -- (-1.743,-1.2949);
\node[solpoint] at (0,0) {};
\node[olab,anchor=north,xshift=14pt,yshift=-1.85pt] at (0,0) {$O$};
\end{tikzpicture}
  \caption{Case (12), $l>a$: the complementary positive construction, with the corresponding change in the conic relations}\label{fig:case12-gt}
  \end{figure}%

  \begin{figure}[htbp]
  \centering
\begin{tikzpicture}[x=1.806cm,y=1.806cm]
\clip (-3.3,-2.2) rectangle (1.25,2.45);
\draw[localaxis] (0,0) -- (-1.7133,0);
\draw[localaxis] (0,-2.2) -- (0,2.45);
\draw[unchosen,domain=-1.45:3.5,samples=161,smooth,variable=\t] plot ({-(\t)},{\t*(\t-1.3)/1});
\draw[chosen,domain=-2.5273:2.5273,samples=141,smooth,variable=\t] plot ({-(0.9+0.4*cosh(\t))},{0.4*sinh(\t)});
\draw[chosen,domain=-2.0015:1.335,samples=161,smooth,variable=\t] plot ({-(exp(\t))},{1-0.5/exp(\t)});
\draw[aseg,line width=4.8pt,yshift=-2.65pt] (0,0) -- (-1.3,0);
\draw[bseg] (0,0) -- (0,1);
\draw[lseg,line width=2pt,yshift=-1.25pt] (0,0) -- (-0.5,0);
\draw[xseg,yshift=1.85pt] (0,0) -- (-1.7133,0);
\draw[yseg] (-1.7133,0) -- (-1.7133,0.7082);
\node[solpoint] at (0,0) {};
\node[olab,anchor=north,xshift=7pt,yshift=-1.85pt] at (0,0) {$O$};
\end{tikzpicture}
  \caption{Case (13), $l<a$: Khayyam's chosen conics and the remaining parabola. A second solution lies on the other branch of the diameter-based hyperbola, which Khayyam does not discuss; Woepcke supplies that branch as a dashed curve, marking it as his own addition, whereas Rashed and Vahabzadeh draw it in the same manner as the branch used in the construction}\label{fig:case13-lt}
  \end{figure}%

  \begin{figure}[htbp]
  \centering
\begin{tikzpicture}[x=3.360cm,y=3.360cm]
\clip (-1.6,-0.9) rectangle (0.6,1.6);
\draw[localaxis] (0,0) -- (-0.7743,0);
\draw[localaxis] (0,-0.9) -- (0,1.6);
\draw[unchosen,domain=-0.8:1.8,samples=161,smooth,variable=\t] plot ({-(\t)},{\t*(\t-0.15)/1});
\draw[chosen,domain=-2.9576:3.0798,samples=141,smooth,variable=\t] plot ({-(0.275+0.125*cosh(\t))},{0.125*sinh(\t)});
\draw[chosen,domain=-1.7918:0.7419,samples=161,smooth,variable=\t] plot ({-(exp(\t))},{1-0.4/exp(\t)});
\draw[lseg,line width=4.8pt,yshift=-2.65pt] (0,0) -- (-0.4,0);
\draw[bseg] (0,0) -- (0,1);
\draw[aseg,line width=2pt,yshift=-1.25pt] (0,0) -- (-0.15,0);
\draw[xseg,yshift=1.85pt] (0,0) -- (-0.7743,0);
\draw[yseg] (-0.7743,0) -- (-0.7743,0.4834);
\node[solpoint] at (0,0) {};
\node[olab,anchor=north,xshift=7pt,yshift=-1.85pt] at (0,0) {$O$};
\end{tikzpicture}
  \caption{Case (13), $l>a$: the same conic relations, with the placement changed by the order of $l$ and $a$. Both conics pass through the vertex of the diameter-based hyperbola and, as Khayyam notes, may meet again by contact or at two points of intersection; the figure shows the latter configuration, with the solution marked at one of the two intersections}\label{fig:case13-gt}
  \end{figure}%

This economy should not be confused with ignorance of alternative conic methods. Khayyam knew a broader mathematical practice in which problems now expressible as cubic equations were treated through different choices of conics. One example is the construction of the regular heptagon, a problem whose Greek and Arabic conic constructions have been documented in detail by Hogendijk \citep{Hogendijk1984}; Khayyam himself discusses the problem in his treatise \citep[p.~174]{RashedVahabzadeh2000}.

This suggests a geometrical explanation for Khayyam's choices. The third conic is not unavailable: it is algebraically available and geometrically constructible. Rather, Khayyam consistently chooses the pair that can be assembled most economically from the local geometrical data already in place. The remaining conic requires an additional geometrical construction before its cross is known in position in relation to the segments already in place.

We can now return to Rashed and Vahabzadeh's appeal to the law of homogeneity. Their explanation does not account for the pattern just described. Homogeneity cannot be the deciding principle: the remaining relation generated by the same proportional scheme is itself homogeneous in Khayyam's sense and defines a legitimate conic relation. What distinguishes Khayyam's chosen pair is geometrical economy. The two conics can be assembled from magnitudes already constructed and known in position, whereas bringing the remaining conic into the same argument requires a further construction to make its cross known in position. The reason for Khayyam's choices is therefore to be sought in the geometry of his procedure, not in homogeneity.

\section{Working inductively within a local framework}

Khayyam himself organizes the cubic equations according to the number of terms and their distribution between the two sides, arriving at thirteen species. For each species he determines a suitable pair of conics, proceeding, in his own words, ``by means of a modicum of induction together with an easy deduction'' \citep[p.~145]{RashedVahabzadeh2000}. To get a sense of what such an inductive process involves within a local geometrical framework, it is useful to reorganize the thirteen species according to the pair of conics Khayyam chooses:

\begin{enumerate}
\item Circle + parabola: Case~(1).
\item Parabola + diameter-based hyperbola: Cases~(2) and (3).
\item Parabola + asymptotic hyperbola: Cases~(4), (5), and (6).
\item Circle + asymptotic hyperbola: Cases~(7), (9), and (12).
\item Diameter-based hyperbola + asymptotic hyperbola: Cases~(8), (10), (11), and (13).
\end{enumerate}

The purpose is emphatically not to reconstruct Khayyam's thought process. Rather, by working through one of these families, we can experience some of the geometrical decisions that must be made when the conics are constructed locally, and see what disappears when the same relations are transferred to a global coordinate framework.

Consider the family formed by Cases~(7), (9), and (12). In all three cases the first conic is a circle. In Cases~(7) and (9) its diameter has length $l+a$. Case~(12) is more revealing: Khayyam treats separately the conditions $l<a$, $l=a$, and $l>a$. In the two nontrivial constructions the diameter is respectively $a-l$ and $l-a$; when $l=a$, he identifies $x=a$ directly. In each construction we must also decide where to place the point $O$, from which the sought red segment $x$ is to be measured.

For Case~(7),
\[
y^2=(l-x)(x+a).
\]
Thus the dashed segment $y$ is the geometric mean of the two parts $l-x$ and $x+a$, exactly as can be read from Fig.~\ref{fig:mean-case7}. This way of putting the guiding conic first is close to Khayyam's own manner of proceeding. He first constructs the semicircle and names on it the point that will later become the intersection with the second conic. In this case he constructs the semicircle $DGC$ on diameter $DC$ and names $G$ as a point on it. At this stage $G$ is simply a point on the first conic; only after the second conic has been constructed does it become the common point required by the solution. In this sense, the point that will eventually be the point of intersection is already present in the first construction before there is a second curve with which to intersect.

  \begin{figure}[htbp]
  \centering
\begin{tikzpicture}[x=3.097cm,y=3.097cm]
\clip (-1.95,-0.55) rectangle (1.15,1.58);
\draw[localaxis] (0.8,0) -- (-1.6,0);
\draw[chosen,domain=0:180,samples=91,smooth,variable=\t] plot ({-(0.4+1.2*cos(\t))},{1.2*sin(\t)});
\draw[lseg,line width=3.2pt,yshift=-1.85pt] (0,0) -- (-1.6,0);
\draw[aseg,line width=3.2pt,yshift=-1.85pt] (0,0) -- (0.8,0);
\draw[xseg,yshift=1.85pt] (0,0) -- (-0.7443,0);
\draw[yseg] (-0.7443,0) -- (-0.7443,1.1495);
\node[solpoint] at (0,0) {};
\node[olab,anchor=north,xshift=7pt,yshift=-5.05pt] at (0,0) {$O$};
\end{tikzpicture}
  \caption{Case (7): the first conic alone. The dashed segment is the geometric mean of the two portions $l-x$ and $x+a$ into which its foot divides the diameter, and its upper end is the point named before the second conic exists}\label{fig:mean-case7}
  \end{figure}%

For Case~(9),
\[
y^2=(x+l)(a-x),
\]
so the dashed segment $y$ is the geometric mean of the two lengths $x+l$ and $a-x$ (Fig.~\ref{fig:mean-case9}). Khayyam's own diagram places the segment $a$ on the left and $l$ on the right (Fig.~\ref{fig:case9}). What matters here is not a fixed left-to-right convention but the way in which the two positive factors of the geometric mean determine the two portions of the diameter; Fig.~\ref{fig:mean-case9} therefore deliberately keeps the order of Case~(7).

  \begin{figure}[htbp]
  \centering
\begin{tikzpicture}[x=1.433cm,y=1.433cm]
\clip (-2.35,-0.55) rectangle (4.35,3.4);
\draw[localaxis] (-2,0) -- (4,0);
\draw[chosen,domain=0:180,samples=91,smooth,variable=\t] plot ({1+3*cos(\t)},{3*sin(\t)});
\draw[aseg,line width=3.2pt,yshift=-1.85pt] (0,0) -- (4,0);
\draw[lseg,line width=3.2pt,yshift=-1.85pt] (0,0) -- (-2,0);
\draw[xseg,yshift=1.85pt] (0,0) -- (1,0);
\draw[yseg] (1,0) -- (1,3);
\node[solpoint] at (0,0) {};
\node[olab,anchor=north,xshift=-7pt,yshift=-5.05pt] at (0,0) {$O$};
\end{tikzpicture}
  \caption{Case (9): the first conic alone, drawn with the order of the given segments reversed relative to Khayyam's own diagram (cf.\ Fig.~\ref{fig:case9})}\label{fig:mean-case9}
  \end{figure}%

Case~(12) makes that dependence on positive order explicit. When $l<a$,
\[
y^2=(x-l)(a-x),
\]
so the circle has positive diameter $a-l$ (Fig.~\ref{fig:mean-case12-lt}); when $l>a$,
\[
y^2=(l-x)(x-a),
\]
the corresponding positive diameter is $l-a$ (Fig.~\ref{fig:mean-case12}). Across these figures, the same geometrical principle can be checked directly: the dashed perpendicular is the geometric mean of the two portions into which its foot divides the diameter, and those portions have precisely the lengths appearing in the corresponding relation. It is worth noticing the literal movement of the eye required to verify this: one follows the coloured diameter from one end, through the foot of the dashed perpendicular, to the other end, repeatedly checking that the required segments occupy the positions prescribed by the relation.

  \begin{figure}[htbp]
  \centering
\begin{tikzpicture}[x=1.939cm,y=1.939cm]
\clip (-4.35,-0.55) rectangle (0.6,1.9);
\draw[localaxis] (0,0) -- (-4,0);
\draw[chosen,domain=0:180,samples=91,smooth,variable=\t] plot ({-(2.5+1.5*cos(\t))},{1.5*sin(\t)});
\draw[aseg,line width=4.8pt,yshift=-2.65pt] (0,0) -- (-4,0);
\draw[lseg,line width=2pt,yshift=-1.25pt] (0,0) -- (-1,0);
\draw[xseg,yshift=1.85pt] (0,0) -- (-3.8063,0);
\draw[yseg] (-3.8063,0) -- (-3.8063,0.7373);
\node[solpoint] at (0,0) {};
\node[olab,anchor=north,xshift=7pt,yshift=-1.85pt] at (0,0) {$O$};
\end{tikzpicture}
  \caption{Case (12), $l<a$: the first conic alone}\label{fig:mean-case12-lt}
  \end{figure}%

  \begin{figure}[htbp]
  \centering
\begin{tikzpicture}[x=1.939cm,y=1.939cm]
\clip (-4.35,-0.55) rectangle (0.6,1.9);
\draw[localaxis] (0,0) -- (-4,0);
\draw[chosen,domain=0:180,samples=91,smooth,variable=\t] plot ({-(2.5+1.5*cos(\t))},{1.5*sin(\t)});
\draw[lseg,line width=4.8pt,yshift=-2.65pt] (0,0) -- (-4,0);
\draw[aseg,line width=2pt,yshift=-1.25pt] (0,0) -- (-1,0);
\draw[xseg,yshift=1.85pt] (0,0) -- (-1.743,0);
\draw[yseg] (-1.743,0) -- (-1.743,1.2949);
\node[solpoint] at (0,0) {};
\node[olab,anchor=north,xshift=7pt,yshift=-1.85pt] at (0,0) {$O$};
\end{tikzpicture}
  \caption{Case (12), $l>a$: the first conic alone---essentially the preceding figure with the roles of $l$ and $a$ interchanged. Khayyam's own configuration is its vertical mirror (cf.\ Fig.~\ref{fig:case12-gt})}\label{fig:mean-case12}
  \end{figure}%

We then have to place the second conic, an asymptotic hyperbola. Here the problem is no longer merely to recognize the appropriate species of conic. Its asymptotes must be situated in relation to the circle and to the segments already known in position. In Case~(7) the second relation is
\[
x(y+b)=bl.
\]
The expression $y+b$ tells us that the distance read as $y$ must be completed by the already constructed blue segment $b$ to reach the relevant asymptote (Fig.~\ref{fig:fit-case7}). The line through the far end of the blue segment, parallel to the corresponding baseline, is therefore part of the hyperbola's local cross even though, as explained above, we did not draw it in Fig.~\ref{fig:case7}.

  \begin{figure}[htbp]
  \centering
\begin{tikzpicture}[x=1.778cm,y=1.778cm]
\clip (-2.8,-1.45) rectangle (2.6,3.1);
\draw[localaxis] (0,-1) -- (-3,-1);
\draw[localaxis] (0.8,0) -- (-1.6,0);
\draw[localaxis] (0,-1) -- (0,3.1);
\draw[chosen,domain=0:180,samples=91,smooth,variable=\t] plot ({-(0.4+1.2*cos(\t))},{1.2*sin(\t)});
\draw[chosen,domain=-1.0561:1.1939,samples=161,smooth,variable=\t] plot ({-(exp(\t))},{-1+1.6/exp(\t)});
\draw[lseg,line width=3.2pt,yshift=-1.85pt] (0,0) -- (-1.6,0);
\draw[aseg,line width=3.2pt,yshift=-1.85pt] (0,0) -- (0.8,0);
\draw[bseg] (0,0) -- (0,-1);
\draw[xseg,yshift=1.85pt] (0,0) -- (-0.7443,0);
\draw[yseg] (-0.7443,0) -- (-0.7443,1.1495);
\node[solpoint] at (0,0) {};
\node[olab,anchor=north,xshift=7pt,yshift=-5.05pt] at (0,0) {$O$};
\end{tikzpicture}
  \caption{Case (7): fitting the asymptotic hyperbola. The clue $y+b$ places the horizontal asymptote at the far end of the blue segment, below the diameter}\label{fig:fit-case7}
  \end{figure}%

Case~(9) has instead
\[
x(y-b)=bl,
\]
so the same displacement by $b$ is read in the other sense (Fig.~\ref{fig:fit-case9}). Case~(12) changes with the parameter order: for $l<a$ the relation is $x(b-y)=bl$ (Fig.~\ref{fig:fit-case12-lt}), while for $l>a$ it is $x(y+b)=bl$ (Fig.~\ref{fig:fit-case12}). The two Case~(12) figures therefore show not merely a relocation of unchanged equations, but two positive geometrical realizations of the same cubic.

  \begin{figure}[htbp]
  \centering
\begin{tikzpicture}[x=1.324cm,y=1.324cm]
\clip (-2.35,-0.6) rectangle (4.9,5);
\draw[localaxis] (0,1) -- (5.1,1);
\draw[localaxis] (-2,0) -- (4,0);
\draw[localaxis] (0,0) -- (0,5);
\draw[chosen,domain=0:180,samples=91,smooth,variable=\t] plot ({1+3*cos(\t)},{3*sin(\t)});
\draw[chosen,domain=-0.8109:1.6864,samples=161,smooth,variable=\t] plot ({exp(\t)},{1+2/exp(\t)});
\draw[aseg,line width=3.2pt,yshift=-1.85pt] (0,0) -- (4,0);
\draw[lseg,line width=3.2pt,yshift=-1.85pt] (0,0) -- (-2,0);
\draw[bseg] (0,0) -- (0,1);
\draw[xseg,yshift=1.85pt] (0,0) -- (1,0);
\draw[yseg] (1,0) -- (1,3);
\node[solpoint] at (0,0) {};
\node[olab,anchor=north,xshift=-7pt,yshift=-5.05pt] at (0,0) {$O$};
\end{tikzpicture}
  \caption{Case (9): the clue $y-b$ places the horizontal asymptote above the diameter, at the far end of the blue segment}\label{fig:fit-case9}
  \end{figure}%

  \begin{figure}[htbp]
  \centering
\begin{tikzpicture}[x=1.939cm,y=1.939cm]
\clip (-4.35,-0.95) rectangle (0.6,2);
\draw[localaxis] (0,1) -- (-4.55,1);
\draw[localaxis] (0,0) -- (-4,0);
\draw[localaxis] (0,-0.95) -- (0,1);
\draw[chosen,domain=0:180,samples=91,smooth,variable=\t] plot ({-(2.5+1.5*cos(\t))},{1.5*sin(\t)});
\draw[chosen,domain=-0.8961:1.579,samples=161,smooth,variable=\t] plot ({-(exp(\t))},{1-1/exp(\t)});
\draw[aseg,line width=4.8pt,yshift=-2.65pt] (0,0) -- (-4,0);
\draw[bseg] (0,0) -- (0,1);
\draw[lseg,line width=2pt,yshift=-1.25pt] (0,0) -- (-1,0);
\draw[xseg,yshift=1.85pt] (0,0) -- (-3.8063,0);
\draw[yseg] (-3.8063,0) -- (-3.8063,0.7373);
\node[solpoint] at (0,0) {};
\node[olab,anchor=north,xshift=7pt,yshift=-1.85pt] at (0,0) {$O$};
\end{tikzpicture}
  \caption{Case (12), $l<a$: the clue $b-y$ places the asymptote above the diameter, in the configuration of Fig.~\ref{fig:mean-case12-lt}}\label{fig:fit-case12-lt}
  \end{figure}%

  \begin{figure}[htbp]
  \centering
\begin{tikzpicture}[x=1.412cm,y=1.412cm]
\clip (-5.6,-1.5) rectangle (1.2,4.2);
\draw[localaxis] (0,-1) -- (-5.8,-1);
\draw[localaxis] (0,0) -- (-4,0);
\draw[localaxis] (0,-1) -- (0,4.2);
\draw[chosen,domain=0:180,samples=91,smooth,variable=\t] plot ({-(2.5+1.5*cos(\t))},{1.5*sin(\t)});
\draw[chosen,domain=-0.3542:1.8083,samples=161,smooth,variable=\t] plot ({-(exp(\t))},{-1+4/exp(\t)});
\draw[lseg,line width=4.8pt,yshift=-2.65pt] (0,0) -- (-4,0);
\draw[bseg] (0,0) -- (0,-1);
\draw[aseg,line width=2pt,yshift=-1.25pt] (0,0) -- (-1,0);
\draw[xseg,yshift=1.85pt] (0,0) -- (-1.743,0);
\draw[yseg] (-1.743,0) -- (-1.743,1.2949);
\node[solpoint] at (0,0) {};
\node[olab,anchor=north,xshift=7pt,yshift=-1.85pt] at (0,0) {$O$};
\end{tikzpicture}
  \caption{Case (12), $l>a$: the clue $y+b$ again places the asymptote below the diameter, in the configuration of Fig.~\ref{fig:mean-case12}}\label{fig:fit-case12}
  \end{figure}%

All of this disappears when one adopts a global Cartesian framework. There we begin by fixing the reference plane. Signed quantities allow the equations themselves to encode on which side of an axis a point or curve lies; translations and reflections are consequently absorbed into the notation. The order of construction, the repeated checking that the relevant segments correspond, and the case-by-case decisions about the placement and orientation of the conics cease to appear as mathematical tasks. The graphs arrive already fitted into a common plane.

\section{Khayyam's programme}

Khayyam's geometry has often been read through the lens of analytic geometry. In doing so, we have inadvertently imposed upon it a global coordinate framework foreign to his method. As a result, his choices of conics have appeared either obvious or unmotivated, and we have lacked the conceptual framework needed to ask why he made them. Yet that framework is present from the outset in Khayyam's own account, explicitly so when he says that the reader needs only Euclid's \emph{Elements} and Apollonius' \emph{Conics}. His objects are geometrical objects, not analytic objects defined in a global coordinate system. Only within this distinction does the central question arise---why these conics rather than others?---and, as we have seen, the answer lies there as well.

In Khayyam's geometry, coordinates are constructed locally, tied to particular diameters, tangents, or asymptotes. Equality remains transitive algebraically: if $a=c$ and $b=c$, then $a=b$. The corresponding local geometrical constructions, however, need not be mutually compatible. The same symbols may represent equal lengths without those lengths already occupying compatible positions in a single figure. What a Cartesian representation achieves automatically must therefore be achieved geometrically in Khayyam's constructions. This is the work hidden behind what Khayyam presents as no more than a little induction followed by an easy deduction. Even where the algebraic relations are mutually compatible, the corresponding conics still have to be placed so that their local measurements become compatible.

This distinction bears directly on a broader historical narrative. Khayyam is often presented as a precursor of Fermat and Descartes. Rashed goes further: in \emph{Classical Mathematics from al-Khwarizmi to Descartes}, the discussion of seventeenth-century geometry appears under the heading ``The geometrical theory of algebraic equations: the completion of al-Khayyam's program'' \citep{Rashed2014}. Yet in matters of coordinate geometry Khayyam stands on fundamentally different ground. His coordinates are local; analytic geometry is global. In his work, geometry and algebra cooperate, but they do not yet merge into a single system in which an algebraic equation determines a curve already situated in a common plane.

That synthesis belongs to the seventeenth century. Fermat could write:

\begin{quote}
Whenever two unknown quantities are found in final equality, there results a locus in place, and the endpoint of one of these describes a straight line or a curve.
\end{quote}

Described by Boyer as ``one of the most significant statements in the history of mathematics'' \citep[p.~75]{Boyer1959}, this marks a fundamental change. Algebra now determines geometry within a fixed global framework: once the variables have been assigned to the axes, an equation determines a locus in the same plane as every other equation expressed in those variables. Translations, changes of orientation, and the reconciliation of local measuring systems become operations internal to the coordinate formalism itself.

A useful historical parallel can be found in Michael N. Fried's discussion of Book IV of Apollonius' \emph{Conics}. Heath had dismissed Book IV as ``on the whole dull''; Fried argues that this dismissal comes from ``deeply rooted prejudices''---in particular, from judging ancient mathematics according to the degree to which it corresponds to modern mathematics when rewritten as a modern mathematical text \citep[p.~267]{Fried2013}. Approached on its own terms, Fried observes, Book IV instead ``raises questions as to Apollonius's basic understanding of how conic sections may be present and related to one another in a single plane'' \citep[p.~268]{Fried2013}.

The parallel with Khayyam is striking. From a Cartesian perspective, much of the problem addressed in Apollonius' Book IV is absorbed into the global framework: conics entered into the same coordinate plane are already situated relative to one another, and their possible intersections become questions about simultaneous equations. What was a geometrical problem of how conics may be present and related in one plane can consequently appear ``dull.'' But this is precisely the kind of problem that remains alive in Khayyam's work. His conics are not born into a common plane. They must be constructed, placed, and related to one another before their intersection can carry the required algebraic meaning.

Seen exclusively through a Cartesian lens, the most conspicuous achievement of Khayyam's cubic theory is therefore his reduction of cubic equations to intersections of conic sections. Yet the use of intersections of conics to solve problems did not originate with him. It had precedents in Greek mathematics, most famously in the duplication of the cube \citep{Heath1921}, and in constructions such as the regular heptagon \citep{Hogendijk1984}. Khayyam's achievement cannot therefore be located simply in the idea of using two conics and taking their intersection.

Once the initial algebraization is placed within a global Cartesian framework, as is often done for educational purposes (for one example among many, see \citealp{KentSherman2015}), much of what follows in Khayyam's treatise is absorbed into the coordinate system itself. The equations place and orient the curves automatically; the geometrical decisions that Khayyam had to make case by case consequently disappear from view.

The same distinction qualifies statements such as Rashed and Vahabzadeh's observation that an equation ``can in certain cases have two negative roots'' \citep[p.~43]{RashedVahabzadeh2000}. As a statement about the corresponding modern algebraic equation, this is unproblematic. It cannot, however, be read as such from Khayyam's diagram. His construction deals with positive magnitudes, and the diagram itself does not establish a signed direction. To distinguish a positive root from a negative one geometrically, a direction---and therefore a convention of sign---must first be imposed.

\begin{figure}[htbp]
\centering
\begin{minipage}[t]{0.538\linewidth}
  \centering
  \includegraphics[width=\linewidth]{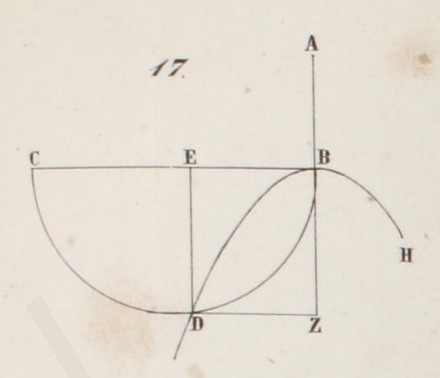}
  \caption{Khayyam's diagram for Case~(1) in Woepcke's edition}
  \label{fig:woepcke-case1}
  {\footnotesize\emph{Source}: \citealp{Woepcke1851}, Plate~II, Fig.~17 (detail)\par}
\end{minipage}\hfill
\begin{minipage}[t]{0.411\linewidth}
  \centering
  \includegraphics[width=\linewidth]{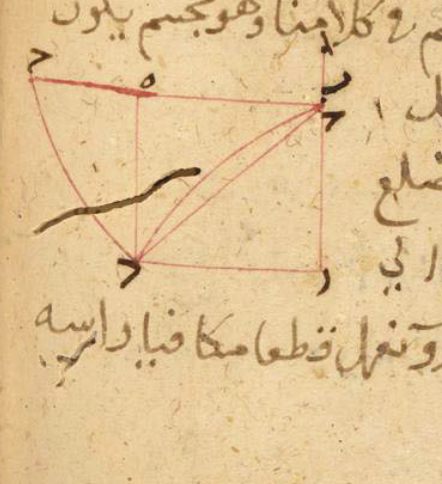}
  \caption{The same diagram in a thirteenth-century manuscript}
  \label{fig:columbia-case1}
  {\footnotesize\emph{Source}: Smith Oriental MS~45, Rare Book \& Manuscript Library, Columbia University\par}
\end{minipage}
\end{figure}

The urge to connect Khayyam to analytic geometry shows itself even in the diagrams Rashed and Vahabzadeh redrew for the treatise. Throughout, the two perpendicular lines that carry the data of a construction are treated as axes, extended beyond what the configuration of the case at hand requires; their diagram for Case~(1) (Fig.~\ref{fig:rashed-reading}) is just one example. There, the line containing the diameter $BC$ of the semicircle is prolonged beyond the semicircle, giving it the visual appearance of an axis. They then make this reading explicit: ``Let axes $(AB,BC)$ be considered as axes $(Oy,Ox)$'' \citep[p.~37]{RashedVahabzadeh2000}. In Woepcke's reproduction of the same construction (Fig.~\ref{fig:woepcke-case1}), by contrast, $BC$ remains simply the diameter of the semicircle; the same is true in the Columbia manuscript (Fig.~\ref{fig:columbia-case1}). The detail is small but systematic, and revealing: lines belonging to Khayyam's geometrical constructions are redrawn as if they were already part of a global coordinate system.

To describe analytic geometry as ``the completion of al-Khayyam's program'' therefore risks concealing as much as it reveals. It makes Khayyam appear as an incomplete analytic geometer, while overlooking the coherence of the geometrical programme he actually pursued. More importantly, it obscures what had to change before analytic geometry could emerge. The decisive transformation was not simply the replacement of verbal algebra by symbolic algebra, nor merely the association of equations with conics. It was the passage from a geometry in which each conic carried its own local frame to one in which all curves could inhabit a single coordinate plane and algebra itself could determine their position within it. Reading Khayyam on his own terms reveals both the richness of his local geometrical constructions and the fundamental change introduced by Fermat and Descartes: the passage to a global coordinate framework.

\section*{Statements and Declarations}

\textbf{Funding} No funding was received for this work.

\textbf{Competing interests} The author has no relevant financial or non-financial interests to disclose.

\textbf{Data availability} No datasets were generated or analysed during the current study. All sources discussed are cited in the references; the figures of Khayyam's constructions are generated from the equations given in the text, and their LaTeX/TikZ source forms part of the manuscript files.

\backmatter
\FloatBarrier
\nocite{Woepcke1851,ColumbiaMS45}
\bibliography{Khayyam_Choices_of_Conics_references}

\end{document}